\documentclass[12pt]{article}%
\usepackage{amsmath,graphicx,amsfonts,amssymb,amsthm,pictexwd}
\providecommand{\U}[1]{\protect\rule{.1in}{.1in}}
\def\theenumi{\arabic{enumi}}

\def\theenumii{\alph{enumii}}
\def\p@enumii{\theenumi.}

\def\theenumiii{\arabic{enumiii}}
\def\p@enumiii{(\theenumi)(\theenumii)}

\def\p@enumiv{\p@enumiii.\theenumiii}
\newcommand{\C}{\mathbb{C}}
\newcommand{\R}{\mathbb{R}}
\newcommand{\V}{\varphi}
\newcommand{\B}{B_n}
\newcommand{\Pn}{P_n}
\newcommand{\CP}{\mathbb{C}P_n}
\newtheorem{theorem}{Theorem}[section]
\newtheorem{corollary}[theorem]{Corollary}
\newtheorem{lemma}[theorem]{Lemma}
\newtheorem{proposition}[theorem]{Proposition}
\theoremstyle{definition}
\newtheorem{definition}[theorem]{Definition}

\newtheorem*{claim}{Claim}
\begin{document}
\title{The Alexander and Jones Polynomials Through Representations of Rook Algebras}

%\email{bigelow@math.ucsb.edu}
%\email{eramos@cmu.edu}
%\email{reyi@ic.sunysb.edu}
\author{Stephen Bigelow\footnote{University of California Santa Barbara, bigelow@math.ucsb.edu}, Eric Ramos\footnote{Funded by NSF grant DMS - 0852065. Carnegie Mellon University, eramos@cmu.edu}, Ren Yi\footnote{Funded by NSF grant DMS - 0852065. Stony Brook University, reyi@ic.sunysb.edu}}
\maketitle
\begin{abstract}

In the 1920's Artin defined the braid group, $B_n$, in an attempt to understand knots in a more algebraic setting. A braid is a certain arrangement of strings in three-dimensional space. It is a celebrated theorem of Alexander that every knot is obtainable from a braid by identifying the endpoints of each string. Because of this correspondence, the Jones and Alexander polynomials, two of the most important knot invariants, can be described completely using the braid group. There has been a recent growth of interest in other diagrammatic algebras, whose elements have a similar topological flavor to the braid group. These have wide ranging applications in areas including representation theory and quantum computation. We consider representations of the braid group when passed through another diagrammatic algebra, the planar rook algebra. By studying traces of these matrices, we recover both the Jones and Alexander polynomials.
\end{abstract}

\section{Introduction}

The Artin braid groups, $B_n$, are the groups generated by $\sigma_1,\sigma_2,\ldots, \sigma_{n-1}$ satisfying the following relations:

\begin{enumerate}
\item $\sigma_i\sigma_j = \sigma_j \sigma_i$ if $|i-j| > 1$
\item $\sigma_i\sigma_j\sigma_i = \sigma_j\sigma_i\sigma_j$ if $|i-j| = 1$
\end{enumerate}

For our purposes, one will find the abstract presentation insufficient. Instead we will work with the following definition found, for example, in \cite{KT}.

\begin{definition}
A geometric braid on $n \geq 1$ strings is a set $b \subset \R^2 \times I$ formed by $n$ disjoint topological intervals called the strings of $b$ such that the projection $\R^2 \times I \rightarrow I$ is increasing and maps each string homeomorphically onto $I$ with 
\begin{eqnarray*}
b\cap(\R^2 \times \{0\}) = \{ (1,0,0),(2,0,0),\ldots, (n,0,0)\} \\
b\cap(\R^2 \times \{1\}) = \{ (1,0,1),(2,0,1),\ldots, (n,0,1)\} 
\end{eqnarray*}
\end{definition}

The operation in this new interpretation is defined by stacking the two braids. An example of this is seen below:\\

\begin{figure}[h]
\includegraphics{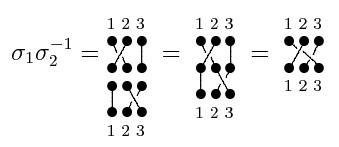}
\end{figure}

These groups are the focus of much study due to their wide ranging applications. One particular area of importance is the representation theory of $\B$. For a comprehensive survey of what is known about these representations, see \cite{BUR}, \cite{FOR} and \cite{BIG}. In this paper, we will study braid group representations when factored through the so called planar rook algebra, $\CP$.\\

The rook monoid, $R_n$, as defined in \cite{FHH}, consists of all bipartite graphs on two rows of $n$ vertices, one on top of the other forming the boundary of a rectangle, such that each vertex has degree either zero or one. Furthermore, we label the vertices in our diagram with the numbers 1 through $n$ from left to right. For example we have the following diagram:

{\beginpicture
\setcoordinatesystem units <0.5cm,0.3cm>        
\setplotarea x from 0 to 6, y from -7 to -.5 
\put{$\bullet$} at 1 -4 \put{$\bullet$} at  1 -1 
\put{$\bullet$} at 2 -4 \put{$\bullet$} at  2 -1
\put{$\bullet$} at 3 -4 \put{$\bullet$} at  3 -1 
\put{$1$} at 1 -.2
\put{$2$} at 2 -.2
\put{$3$} at 3 -.2
\put{$1$} at 1 -4.8
\put{$2$} at 2 -4.8
\put{$3$} at 3 -4.8
\plot 1 -4 3 -1 /
\plot 2 -4 1 -1 /

\put {$\in R_3$} at 5 -2.5
\endpicture}
We also have the following definitions which will simplify discussions of rook diagrams:
\begin{definition}
Given any rook diagram, $d \in R_n$ we say:
\begin{enumerate}
\item $d$ has a \textit{vertical line} at $i$ if it contains an edge from the vertex $i$ on top to the same vertex $i$ on bottom.
\item $d$ has a \textit{dead end} at $i$ if vertex $i$ on top (or bottom) is not incident to an edge.
\item $d$ has a \textit{slash} from $i$ to $j$ if vertex $i$ on bottom is adjacent to vertex $j$ on top.
\end{enumerate}
\end{definition}

The product, $d_1d_2$, of two rook diagrams $d_1$ and $d_2$ is obtained by stacking $d_1$ on top of $d_2$ and deleting any edge from one that connects to a dead end from the other. For example,

{\beginpicture
\setcoordinatesystem units <0.5cm,0.3cm>        
\setplotarea x from 0 to 20, y from -10 to 3
\put{$\bullet$} at 1 -4 \put{$\bullet$} at  1 -1 
\put{$\bullet$} at 2 -4 \put{$\bullet$} at  2 -1
\put{$\bullet$} at 1 -7 \put{$\bullet$} at  2 -7
\put{$\bullet$} at 3 -4 \put{$\bullet$} at  3 -1 \put{$\bullet$} at 3 -7 
\plot 1 -4 2 -1 /
\plot 1 -4 3 -7 /
\plot 3 -1 3 -4 /
\plot 8 -3 9 -5 /
\put{$\bullet$} at 7 -3 \put{$\bullet$} at  7 -5
\put{$\bullet$} at 8 -5 \put{$\bullet$} at  8 -3
\put{$\bullet$} at 9 -5 \put{$\bullet$} at  9 -3
\put{$=$} at 5 -4
\endpicture}

With the rook monoid defined, one considers the submonoid of planar rook diagrams, $\Pn$. These are those diagrams that can be drawn in the plane without crossings when their edges are not allowed to leave the rectangle formed by their vertices. Notice that the diagrams being multiplied above are in $P_3$. Finally, one simply says that the planar rook algebra, $\CP$, is the $\C$-algebra generated by $\Pn$ with multiplication defined using the distributive law.\\

What makes $\CP$ such an ideal object of study, aside from its relatively simple presentation, is that its representations have already been completely classified. One will find this classification in \cite{FHH}. This is invaluable for studying representations of $\B$ through $\CP$. We discover through this paper that many of these representations are those of the Hecke Algebra. We define:
\begin{definition}
The Hecke Algebra, $\mathcal{H}_n(q)$, is the quotient algebra of $\C \B$ by the subalgebra generated by the following relations:
\[
(\sigma_i - 1)(\sigma_i + q) = 0,  \quad i\in \{1,\ldots,n\}.
\]
\end{definition}

This algebra has shown itself to be spectacularly useful in various mathematical fields. In \cite{JON2}, Jones illustrated how this algebra and its representations could be used to fully recover the Jones polynomial invariant. To understand this, one must first understand the connection between links and braids.\\

Looking back to the geometric definition of a braid, it may become clear that there exists some relationship between the braid groups and link theory. Indeed it is a celebrated theorem of Alexander that these two concepts are essentially the same. To make this precise we first define for any braid $b\in \B$, the closure of $b$, denoted $\hat{b}$, is the link in $\R^3$ obtained by identifying the top and bottom of $b$. From this definition one can prove the following theorem from \cite{ALE1} and \cite{ALE2}.
\begin{theorem}[Alexander]
Every oriented link is isotopic to a closed braid
\end{theorem}

Unfortunately, the above correspondence is not one to one. Though this may seem like a major setback, there is a way to exactly characterize the correspondence's failure of bijectivity using a method described first by Markov. We define an equivalence relation, $\sim$, on braids as being generated by the following three ``Markov moves":

\begin{enumerate}
\item $ab \sim ba$
\item $\iota(a)\sigma_{n} \sim a$
\item $\iota(a)\sigma_{n}^{-1} \sim a$
\end{enumerate}
where $a,b\in B_n$ and $\iota: B_n\hookrightarrow B_{n+1}$ is the map that adds a vertical string to the end of a braid. Using this one then proves the following theorem.

\begin{theorem} [Markov]
Given two braids, $b_1$ and $b_2$, $\hat{b}_1 = \hat{b}_2$ if and only if $b_1 \sim b_2$
\end{theorem}

In this paper we will be expanding on the ideas found in \cite{BIG3} by creating and using Markov traces to emulate the results of Jones in \cite{JON1}. We find that both the Jones and Alexander Polynomials arise through traces of the Planar Rook Algebra.

\section{Finding the homomorphism}

Before we begin, it is important that we impose an ordering on $P_2$ as this will allow for much simpler notation later on.

\begin{definition}
We enumerate the elements of $P_2$ in the following fashion:
\end{definition}

{\beginpicture
\setcoordinatesystem units <0.5cm,0.3cm>        
\setplotarea x from 0 to 30, y from -8 to 2
\put{$d_1 = $} at 0 -2.5
\put{$\bullet$} at 1.5 -4 \put{$\bullet$} at  1.5 -1 
\put{$\bullet$} at 2.5 -4 \put{$\bullet$} at  2.5 -1

\put{$d_2 = $} at 5 -2.5
\put{$\bullet$} at 6.5 -4 \put{$\bullet$} at  6.5 -1 
\put{$\bullet$} at 7.5 -4 \put{$\bullet$} at  7.5 -1
\plot 6.5 -1 6.5 -4 /

\put{$d_3 = $} at 10 -2.5
\put{$\bullet$} at 11.5 -4 \put{$\bullet$} at  11.5 -1 
\put{$\bullet$} at 12.5 -4 \put{$\bullet$} at  12.5 -1
\plot 11.5 -4 12.5 -1 /

\put{$d_4 = $} at  15 -2.5
\put{$\bullet$} at 16.5 -4 \put{$\bullet$} at  16.5 -1 
\put{$\bullet$} at 17.5 -4 \put{$\bullet$} at  17.5 -1
\plot 16.5 -1 17.5 -4 /

\put{$d_5 = $} at 20 -2.5
\put{$\bullet$} at 21.5 -4 \put{$\bullet$} at  21.5 -1 
\put{$\bullet$} at 22.5 -4 \put{$\bullet$} at  22.5 -1
\plot 22.5 -4 22.5 -1 /

\put{$d_6 = $} at 25 -2.5
\put{$\bullet$} at 26.5 -4 \put{$\bullet$} at  26.5 -1 
\put{$\bullet$} at 27.5 -4 \put{$\bullet$} at  27.5 -1
\plot 26.5 -4 26.5 -1 /
\plot 27.5 -4 27.5 -1 /
\endpicture}

Next, we present the homomorphism which we will be working with for the remainder of the paper, $\varphi:B_n\rightarrow (\CP)^*$ where $(\CP)^*$ is the group of units of $\CP$. To do this, we introduce the following definition:

\begin{definition}
Given two diagrams, $a,b\in P_n$, we define the tensor product, denoted $a\otimes b$, to be the result of appending of $b$ to the right of $a$.
\end{definition}

With these in mind we begin looking for a homomorphism of the following form:
\[
\V(\sigma_i) = a\cdot d_{1i} + b\cdot d_{2i} + c\cdot d_{3i} + d\cdot d_{4i} + e\cdot d_{5i} + f \cdot d_{6i} 
\]
where 
\[d_{ji} = I^{\otimes i-1} \otimes d_j \otimes I^{\otimes n-i-1}
\]
and $I$ is the identity diagram in $P_1$. This is shown explicitly below:
\begin{figure}[h]
\includegraphics{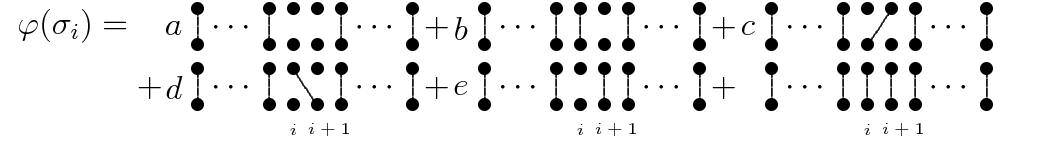}
\end{figure}

It therefore remains to find the conditions on our coefficients which make $\V$ a group homomorphism. In particular we need to make sure the braid relations are satisfied. Furthermore, we need to make sure every element mapped to in $\CP$ is invertible. For this reason we claim $f$ cannot be zero. If it was, our homomorphism would not have an identity diagram as one of its terms. One then notices that it is impossible for the identity to be created by multiplication of other diagrams. Because we require $\V(\sigma_i^{-1}) = \V(\sigma_i)^{-1}$, it follows that $f$ must be non-zero. Therefore, for the purpose of simplicity we will scale $f$ to be 1. After doing this we find the following:

\begin{theorem}\label{homom}
Assuming $a+c+d \neq 1, f=1$ and $cd\neq 0$, any mapping of the above form is a homomorphism if and only if its coefficients are in one of the following five families:
\begin{enumerate}
\item $ b = e = -1$
\item $ a = -c-d$, $b = -1$, $e = -cd$
\item $ a = -c-d$, $b = -cd$, $e = -1$
\item $ a = 1-c-d+ cd$, $b = -cd$, $e =-1 $
\item $ a = 1-c-d+cd$, $b = -1$, $e = -cd$
\end{enumerate}
Furthermore, for each of the above five families one finds the following
\begin{enumerate}
\item
$\varphi(\sigma_i^{-1}) = \left(1 - \frac{1}{d} - \frac{1}{c} + \frac{1}{-1+a+c+d}\right)d_{1i} - d_{2i} + \frac{1}{d}d_{3i} + \frac{1}{c}d_{4i}- d_{5i}+ d_{6i}$
\item
$\varphi(\sigma_i^{-1}) = (-\frac{1}{c} - \frac{1}{d})d_{1i} - \frac{1}{cd}d_{2i} + \frac{1}{d}d_{3i} +\frac{1}{c}d_{4i} - d_{5i}  + d_{6i}$ 
\item 
$\varphi(\sigma_i^{-1}) = (-\frac{1}{c} - \frac{1}{d})d_{1i} - d_{2i} + \frac{1}{d}d_{3i} +\frac{1}{c}d_{4i} -\frac{1}{cd} d_{5i}+ d_{6i}$
\item
$\varphi(\sigma_i^{-1}) = (1 -\frac{1}{c} - \frac{1}{d} + \frac{1}{cd})d_{1i} - d_{2i} + \frac{1}{d}d_{3i} + \frac{1}{c}d_{4i}- \frac{1}{cd}d_{5i} + d_{6i}$
\item
$\varphi(\sigma_i^{-1}) = (1 -\frac{1}{c} - \frac{1}{d} + \frac{1}{cd})d_{1i} -\frac{1}{cd} d_{2i} + \frac{1}{d}d_{3i}  + \frac{1}{c}d_{4i}- d_{5i} + d_{6i}$
\end{enumerate}
\end{theorem}
\begin{proof}
We begin by considering the first braid relation. If two generators, $\sigma_i$ and $\sigma_j$, are such that $|i-j| > 1$, then the diagrams in their image will not overlap in their copies of $d_i$. In particular, in any product of diagrams the $d_i$ will always end up on top, or bottom, of vertical lines. It then follows that this will commute. We find, therefore, the first braid relation will hold regardless of what coefficients we choose.\\

Next we claim it suffices consider the second braid relation for $\V: B_3 \rightarrow (\C P_3)^*$. One will notice that in the product, $\V(\sigma_i)\cdot\V(\sigma_{i+1})\cdot \V(\sigma_i)$, the only parts of the diagrams which actually change are the copies of the $d_j$. Everything else in these diagrams are vertical lines at the ends, which act like identities. From this we see that for higher values of $n$, one will only be adding more vertical lines which do not influence anything.\\

Once these facts have been established, one may proceed by exhaustive case study in proving the first half of the theorem. This can be done by a computer by considering the regular representation of $\CP$ and identifying the matrix $\V(\sigma_1\sigma_{2}\sigma_{1}) - \V(\sigma_2\sigma_{1}\sigma_{2})$ with the zero matrix. This results in 400 polynomials that are required to equal zero. We factored these polynomials, and observed that the expression $c(b+1)(e+1)$ appeared frequently. We therefore solved the system of equations separately in the cases where c, b+1 or e+1 is zero. After finishing this, any redundant or non-invertible solutions were removed and our result follows. \\

Once this has been complete, one simply inverts the matrices found for $\V(\sigma_1)$ to finish the proof.

\end{proof}
For convenience, we identify the above homomorphisms by $\V_i$, where $i=1,...,5$. One may notice the symmetry in the variables $b$ and $e$ as well as $c$ and $d$. Intuitively one may justify this by considering what happens to a diagram when rotated or reflected. To be more precise, we first notice that in each pair of ``dual" families, the only free variables are $c$ and $d$. If we broadcast our choice by saying for $i\in \{2,\ldots,5\}, \V_i = \V_i^{c,d}$ then the following corollary becomes apparent.

\begin{corollary}
Let $\alpha: B_n \rightarrow B_n$ be the automorphism which sends $\sigma_i$ to $\sigma_i^{-1}$
\begin{enumerate}
\item $\V_2^{c,d}(\sigma_i) =  \V_3^{\frac{1}{d}, \frac{1}{c}}(\alpha(\sigma_i))$
\item $\V_5^{c,d}(\sigma_i) = \V_4^{\frac{1}{d},\frac{1}{c}}(\alpha(\sigma_i))$
\end{enumerate}
\end{corollary}

Because of this corollary, we are free to only consider $\V_1$, $\V_2$ and $\V_5$. For the remainder of this paper we will be looking at each one of these three morphisms and discovering representations which factor through them as well as knot invariants. We find that $\V_2$ and $\V_5$ have many very similar, and sometimes exactly the same, properties. We will also find that $\V_1$ is the most different and will therefore be treated last. With all this said we begin with $\V_5$. 

\section{The Hecke algebra and Jones Polynomial}
%%%%%%%%%%%%%%%%%%%%
We begin this section by quickly considering the representations that factor through $\V_5$. Once this is completed we will discuss any knot invariants that can be found using this homomorphism.
\subsection{Representations Through $\V_5$}
We begin our study by showing $\V_5$ satisfies a HOMFLYPT polynomial. Before we begin however, it is important one take note of the form it is presented in. In particular, we show a skein relation in terms of the braids $\sigma_i, \sigma_i^{-1}$ and $1$. It is non-obvious that any pair of links that are related by a crossing change in an arbitrary diagram can be obtained from braid closures in this way. For a proof of this fact one must simply look to the proof of Alexander's theorem itself. This can be found, for example, in \cite{KT}. After one understands this algorithm for converting a braid to a link it becomes clear that what we are working with is indeed equivalent. \\

With all of this said we proceed with the following theorem:
\begin{lemma}\label{skein}
The homomorphism, $\V_5: B_n \rightarrow (\CP)^*$, satisfies the following skein relation for all $\sigma_j$, $j\in \{1,\ldots,n-1\}$ and $x\in{\B}$
\[
\V_5 (x\sigma_j) - cd \cdot \V_5 (x\sigma^{-1}_j) = (1-cd) \cdot \V_5(x).
\] 
In particular these homomorphisms satisfy the following quadratic equation
\[
(\V_5(\sigma_j) - \V_5(1))\cdot(\V_5(\sigma_j) + cd\cdot \V_5(1)) = 0.
\]
\end{lemma}

\begin{proof}
We begin by noticing it suffices to prove for both assertions that,
\[
\V_5 (\sigma_j) - cd \cdot \V_5 (\sigma^{-1}_j) = (1-cd) \cdot \V_5(1).
\] 
We proceed by gathering like terms to find
\[
(1-c-d+cd - cd(\frac{1-c-d+cd}{cd}))d_{1j} + (-1 + \frac{cd}{cd})d_{2j} + (c - \frac{cd}{d})d_{3j} + (d - \frac{cd}{c})d_{4j} + (-cd + cd)d_{5j} + (1-cd)d_{6j}.
\]
Simplifying the above gives our result.
\end{proof}

From this we may immediately categorize the representations of $B_n$ through $\V_5$. To do this one simply combines the quadratic equation mentioned in the statement of Lemma \ref{skein} in conjunction with the definition of the Hecke algebra. To be precise,

\begin{theorem}
All representations of $\B$ through $\V_5$ are representations of $\mathcal{H}_n(cd)$
\end{theorem}

\subsection{Rediscovering the Jones polynomial}
Now that those representations which factor through $\V_5$ have been shown to fall into a well studied class of representations, we turn our attention to the question of knot invariants. All of our results in this topic require the use of trace functions.
\begin{definition}
 A \textit{trace function}, $tr:A\rightarrow F$ is any linear function from an algebra, $A$, to a field, $F$, which satisfies $tr(ab) = tr(ba)$
\end{definition}

It is a simple exercise to show that if $A$ is the algebra of all $n\times n$ of matrices then the only trace functions  are scalar multiples of the classical matrix trace. Furthermore, we say a trace is a \textit{Markov trace} if it is invariant under the Markov moves. Using Markov trace functions, we will discover knot invariants associated with each of our homomorphism families.\\

We begin by defining what we call the bubble trace, $tr^{\beta}_n$.

\begin{definition}
For any $\beta \in \C$, the bubble trace function $tr^{\beta}_n:\CP\rightarrow \C$ is the linear function which acts on diagrams, $d\in\Pn$, by $tr^{\beta}_n(d) = \beta^{k(d)}$ where $k(d)$ is the number of vertical lines in $d$.
\end{definition}

Of course it requires some justification to show that this is a trace. We have the following lemma:
\begin{lemma}\label{vli}
For all diagrams $a,b \in \Pn$, $k(ab) = k(ba)$, where $k(d)$ is the number of vertical lines in $d$.
\end{lemma}
\begin{proof}
One considers how a vertical line can be formed in a product. First, one may have two vertical lines stacked on top of one another. This configuration is obviously preserved if product was done in the opposite order. Secondly, one may have a slash from $i$ to $j$ in $a$ and an slash from $j$ to $i$ in $b$. While the location of the vertical line is not preserved by changing the order of $a$ and $b$, its existence is.
\end{proof}
Now that we see $tr_n^{\beta}$ is a trace, we will attempt to make it a Markov trace. This is achieved below:

\begin{proposition}\label{traces}
Let $tr^{\beta}_n$ be a bubble trace, then the following is a Markov trace:
\[
Tr^5_n(x) = (\sqrt{cd})^{w(x)+n} \cdot tr^{\frac{1+cd}{cd}}_n(\V_5(x))
\]
where $w(x)$ is the exponent sum, or writhe, of $x$.
\end{proposition}
\begin{proof}
We begin by noting the first Markov move is preserved by Lemma \ref{vli}. We then consider $tr^{\beta}_{n+1}(\V(x\sigma_n))$ where $x\in B_{n+1}$ does not contain $\sigma_n$. We may visualize this using Figure 1 below.\\

\begin{figure}
\includegraphics{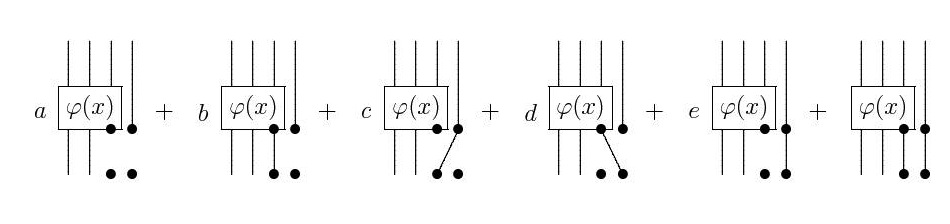}
\caption{}
\end{figure}

From this picture we see that we may gather terms in the following way:
\[
(a+c+d+e\beta)tr^{\beta}_n(\V_5(x)') + (\beta + b)tr^{\beta}_n(\V_5(x))
\]
where $\V(x)'$ is a linear combination of all the diagrams of $\V_5(x)$ but with any vertical lines at $n$, broken. We see that by construction $a+c+d + e\beta = 0$ and so we are left with:
\begin{eqnarray}
tr^{\frac{1+cd}{cd}}_{n+1}(\V_5(x\sigma_n)) = \frac{1}{cd}tr^{\frac{1+cd}{cd}}_n(\V_5(x)) \label{tra}.
\end{eqnarray}
Multiplying \ref{tra} by $(\sqrt{cd})^{w(x\sigma_n)+n+1}$ and using the definition of $Tr_n^5$
\[
Tr^5_{n+1}(x \sigma_n) = (\sqrt{cd})^{w(x)+1+n+1-2}tr^{\frac{1+cd}{cd}}_n(\V_5(x)) = (\sqrt{cd})^{w(x)+n}tr^{\frac{1+cd}{cd}}_n(\V_5(x)) = Tr^5_n(x).
\]
The proof for the last Markov move follows similarly.
\end{proof}

Using this new Markov trace function we will discover that one can indeed recover the Jones polynomial from $\V_5$.  
\begin{theorem}\label{jones}
If $cd \neq -1$ then $Tr^5_n(x) = \frac{1+cd}{\sqrt{cd}} V(\hat{x})$ where $V$ is the Jones polynomial.
\end{theorem}

\begin{proof}
Since $Tr_n^5(x)$ is a Markov trace, we know it is an invariant of the oriented knot, or link, $\hat{x}$. Furthermore, we note the only invariants on knots that satisfy the Jones skein relation are scalar multiples of the Jones polynomial. Therefore we begin by showing that $Tr^5_n$ satisfies the proper skein relation. Once this is finished it will only remain to show that our trace also sends the unknot to $\frac{1+cd}{\sqrt{cd}}$ . One recalls Lemma \ref{skein} and applies $tr^{\frac{1+cd}{cd}}_n$ to both sides to obtain 
\[
tr^{\frac{1+cd}{cd}}_n(\V_5(x\sigma_i)) - cd\cdot tr^{\frac{1+cd}{cd}}_n(\V_5(x\sigma_i^{-1})) = (1-cd)tr^{\frac{1+cd}{cd}}_n(\V_5(x)).
\]
substituting the definition of $Tr^5_n$ given in Proposition \ref{traces} we find,

\begin{description}
\item  $(\sqrt{cd})^{-(n+w(x)+1)}Tr^5_{n}(x\sigma_i) - cd \cdot (\sqrt{cd})^{-(n+w(x)-1)}Tr^5_n(x\sigma^{-1}_i) =(1-cd) \cdot(\sqrt{cd})^{-(n+w(x))}Tr_n^5(x)$
\item $\implies \frac{1}{cd}Tr^5_n(x\sigma_i) - cd\cdot Tr^5_n(x\sigma_i^{-1}) = (\frac{1}{\sqrt{cd}}-\sqrt{cd})\cdot Tr^5_n(x)$
\end{description}

as needed. Furthermore we see in $B_1$,
\[
Tr^5_1(1) = (\sqrt{cd})^{0 + 1}\frac{1+cd}{cd} = \frac{1+cd}{\sqrt{cd}}.
\]
One notices that by our assumptions $cd \neq 0, -1$  and this concludes the proof.
\end{proof}

We may now discuss the second homomorphism family.

\section{The Hecke algebra and the Alexander polynomial}
As stated previously, one will find the results in this section to be very similar to those found above. Before we prove various properties of this homomorphism, we will rescale $\V_2$ for convenience. In particular we scale $f$ to be $\frac{1}{\sqrt{cd}}$ and find:

\begin{eqnarray}
\V_2(\sigma_j) &=&(\frac{-\sqrt{c}}{\sqrt{d}} - \frac{\sqrt{d}}{\sqrt{c}})d_{1j} - \frac{1}{\sqrt{cd}}d_{2j} + \frac{\sqrt{c}}{\sqrt{d}}d_{3j}+\frac{\sqrt{d}}{\sqrt{c}}d_{4j}-\sqrt{cd}\cdot d_{5j} + \frac{1}{\sqrt{cd}}d_{6j} \label{eqn1}.\\
\V_2(\sigma_j^{-1}) &=& (\frac{-\sqrt{c}}{\sqrt{d}} - \frac{\sqrt{d}}{\sqrt{c}})d_{1j} - \frac{1}{\sqrt{cd}}d_{2j} + \frac{\sqrt{c}}{\sqrt{d}}d_{3j}+\frac{\sqrt{d}}{\sqrt{c}}d_{4j}-\sqrt{cd}\cdot d_{5j} + \sqrt{cd}\cdot d_{6j}\label{eqn2}.
\end{eqnarray}

One will notice that all of the coefficients between the two linear combinations above are the same aside from the last one. For the rest of this section we shall refer to our newly scaled homomorphism by $\V_2$.

\subsection{Representations through $\V_2$}

We begin similarly to the $\V_5$ case by showing $\V_2$ satisfies a HOMFLYPT polynomial. In particular we find,

\begin{lemma} \label{skein2}
For each $x\in \B$
\[
\V_2(x\sigma_i) - \V_2(x\sigma_i^{-1}) = (\frac{1}{\sqrt{cd}} - \sqrt{cd})\cdot\V_2(x).
\]
\end{lemma}
\begin{proof}
This follows from equations \ref{eqn1} and \ref{eqn2}.
\end{proof}
Once again the above skein relationship tells us exactly what representations will look like through $\V_2$. We have the following theorem:
\begin{theorem}
All representations of $\B$ through $\V_2$ are representations of $\mathcal{H}_n(1+cd)$
\end{theorem}

\subsection{Rediscovering the Alexander polynomial}

Now that the representations are classified we may begin considering knot invariants. One may expect this to be done using the methods discussed for the previous homomorphism. Unfortunately, this cannot be the case. Going through the same steps above one finds a Markov bubble trace with $\beta = 0$ for $\V_2$. It then can be shown that this trace satisfies the Alexander skein relation. A quick computation reveals that the unknot is sent to zero and thus the work is wasted. We do not give up on our goal and thus define the following trace:

\[ 
tr_n(x) = 
\begin{cases}
1 &\text{ if } x \text{ has exactly 1 vertical line }\\
0 &\text{ otherwise}
\end{cases}.
\]

One may quickly note that once again Lemma \ref{vli} tells us this is indeed a trace. Furthermore, Lemma \ref{skein2} shows that it will satisfy the Alexander skein relation. Because we do not know whether $tr_n$ is a Markov trace, we may not yet conclude it is the Alexander polynomial. The remainder of this subsection will be dedicated to normalizing $tr_n$ to make it a Markov trace. We begin first with the following critical proposition
%%%%%%%%%%%%%%%%%%%%%%%%%%%%%%%%%%%%%%%%%%%%%%
\begin{proposition}\label{vip}
Let $x_n = \sigma_1\cdots\sigma_{n-1} \in \B$. Then the following three statements hold
\begin{enumerate}
\item The sum of all coefficients of diagrams in $\V_2(x_n)$ with no vertical lines is 0.
\item The sum of all coefficients of diagrams in $\V_2(x_n)$ with a vertical line at $n$ and nowhere else is $(-\sqrt{cd})^{n-1}$.
\item $tr_n(\V_2(x_n)) = \left(\frac{-1}{\sqrt{cd}}\right)^{n-1}[1 + cd + (cd)^2 + \ldots + (cd)^{n-1}]$.
\end{enumerate}
\end{proposition}

\begin{proof}
We prove the lemma by induction. \\
If $n=2$ then
\[
\V_2(\sigma_1) = (\frac{-\sqrt{c}}{\sqrt{d}} - \frac{\sqrt{d}}{\sqrt{c}})d_1 - \frac{1}{\sqrt{cd}}d_2 + \frac{\sqrt{c}}{\sqrt{d}}d_3+\frac{\sqrt{d}}{\sqrt{c}}d_4-\sqrt{cd}\cdot d_5 + \frac{1}{\sqrt{cd}}d_6.
\]
One begins by quickly verifying that the coefficients on $d_1$, $d_3$ and $d_4$, the only terms with no vertical lines, sum to zero. Next, one sees that the only diagram with a single vertical line at $n$ is $d_5$ and its coefficient is $-\sqrt{cd}$. For the final claim we see only the fifth and second terms have exactly one vertical line. Summing their two coefficients gives $-\frac{1}{\sqrt{cd}} - \sqrt{cd}$, as needed.\\

Assume that the statement is true for some $n \geq 2$, and consider $x_{n+1}$. We write this as $x_n\sigma_{n}$, where $x_n \in B_{n+1}$, and consider the diagrams found in Figure 1. We see in this case that the first, third and fourth terms $-$ those associated with the empty diagram and the two slashes, respectively $-$ do not gain a vertical line on their far right. We may therefore conclude that these will have no effect in proving claim two. The inductive hypothesis also tells us that any items in these terms with no straight lines will sum to zero, so we need not worry about them for the remainder of claim 1. Finally, when one takes the trace, one finds, similarly to that in the proof of Proposition \ref{traces}, terms will gather in the following way:
\[
((\frac{-\sqrt{c}}{\sqrt{d}} - \frac{\sqrt{d}}{\sqrt{c}}) + \frac{\sqrt{c}}{\sqrt{d}}+ \frac{\sqrt{d}}{\sqrt{c}})tr_n(\V_2(x)').
\]
This sum is zero and thus we are free to ignore these terms for the remainder of the proof.\\

We next turn our attention to the second, fifth, and sixth terms $-$ those associated with the left vertical line, right vertical line, and the identity diagram, respectively. We notice the trace of the sixth term will cause any diagram in $\V_2(x_n)$ with a vertical line to vanish. Furthermore, the terms without a vertical line will sum to zero by the hypothesis. We therefore only need consider the second and fifth terms.\\

Looking then at the fifth term, we see that any terms in $\V_2(x_n)$ with no vertical lines will once again sum to zero. The only terms we consider are those in $\V_2(x_n)$ with exactly one vertical line in the far right, as the break caused by $d_5$ will cause this line to be lost in $\V_2(x_n\sigma_n)$. By our induction hypothesis we know these terms' coefficients sum to $(-\sqrt{cd})^{n-1}$. This completes the proof of claim 2 as the coefficient of the fifth term is $-\sqrt{cd}$. To finish the entire induction one simply notices the second term is exactly $tr(\V(x_n))$ with a coefficient of $\frac{-1}{cd}$. Our inductive hypothesis gives us both our results. \\

\end{proof} 

With this proven we define a new trace in the following way, for any $x\in\B$:

\[
Tr_n^2(x) = \frac{tr_n(\V_2(x))}{tr_n(\V_2(x_n))}.
\]

One will notice that this trace will continue to satisfy the Alexander skein relation.\\

We claim that the above trace is indeed a Markov trace. It is non-obvious how one can prove this directly, however. Instead we will provide a work around using ideas introduced in \cite{BIG2}. For our next statements, we must introduce some notation used in the aforementioned paper.
\begin{definition}
For any partition, $\lambda = (n_1,\ldots,n_k)$, of $n$ we say $\tau_{\lambda}$ is the braid $(\sigma_{1}\cdots\sigma_{n_1-1})\otimes \ldots \otimes (\sigma_{1}\cdots\sigma_{n_k-1})$ where the $i$th term in the product is a braid in $B_{n_i}$.
\end{definition}
To be clear one will find $\tau_{(4,1)}$ below:\\
\begin{center}
\includegraphics{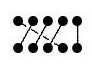}
\end{center}

One may notice that if $\lambda = (n)$ then $\tau_{(n)}$ is exactly the $x_n$ in Proposition \ref{vip}. With this made clear we have another lemma as well as a critical proposition.

\begin{lemma}\label{unl}
if $\lambda \neq (n)$ then $Tr(\tau_{\lambda})=0$.
\end{lemma}

\begin{proof}
We will prove the proposition for the case where $\lambda$ has two parts. One will see that the method used generalizes quite easily. We see that $\tau_{\lambda}=(\sigma_1 \cdots \sigma_i)(\sigma_{i+2} \cdots\sigma_{n})$ for some $i > 0$. Let $\beta_1=\sigma_1\cdots \sigma_i$, $\beta_2=\sigma_{i+2} \cdots\sigma_{n}$. Since $\beta_1$ and $\beta_2$ are disjoint, all of the non vertical portions of diagrams which arise from $\V_2(\beta_1)$ will end up on top of vertical lines in the diagrams of $\V_2(\beta_2)$. From this it follows that when applying $Tr_n^2$ to the product, the surviving diagrams will be a product of a diagram with a single vertical line from one of $\V_2(\beta_1)$ or $\V_2(\beta_2)$, with a diagram with no vertical lines from the other. Once again by Proposition \ref{vip} we know these coefficients sum to zero and we are done.
\end{proof}

\begin{proposition}\label{big}
For all braid words, $w \in \B$, there exists a sequence $w = w_0, w_1,\ldots, w_k$ such that the following hold:
\begin{enumerate}
\item $w_j = x\sigma_i^{\pm1}y$, $w_{j+1} = x\sigma_i^{\mp1}y$ for some $x,y\in \B$ and some $i$.
\item $w_k$ is conjugate to $\tau_{\lambda}$ for some partition $\lambda$
\end{enumerate} 
\end{proposition}
\begin{proof}
the proof can be found in the form of an algorithm discussed in \cite{BIG2}
\end{proof}

With all of these tools we are ready to prove that $Tr_n^2$ is indeed the Alexander polynomial. We have:

\begin{theorem}
If $w\in \B$ then $Tr_n^2(w) = \Delta(\hat{w})$, where $\Delta$ is the Alexander polynomial.
\end{theorem}
\begin{proof}
We prove the theorem by induction on the length of the word $w$. If $|w| = 0$ then $Tr_n^2(\V_2(w)) = 0$, unless $n = 1$ or $2$. This is because in the case of $n = 3$, the extra vertical line makes it so that only terms with no vertical lines will survive the trace. We know these terms sum to zero. For all greater $n$, every term has at least two vertical lines. In the cases of $n = 1$ and $2$, the trace will be 1 as desired by simple calculation.\\

Next assume the statement holds for all lengths up to and including $l$ and let $w$ be an arbitrary word of length $l+1$. To proceed we prove the following claim:
\begin{claim}
All braids in the sequence $w = w_0,\ldots,w_k$ granted by Proposition \ref{big} have $Tr_n^2(w_i) = \Delta(\hat{w}_i)$. 
\end{claim}
\begin{proof}
We prove that $Tr_n^2(w_{k-r}) = \Delta(\hat{w}_{k-r})$ by induction on $r$.\\
if $r=0$ then,
\[
 Tr_n^2(w_{k-r})=Tr_n^2(w_k)=Tr_n^2(\tau_{\lambda})=\Delta(\hat{\tau}_{\lambda})=\Delta(\hat{w}_{k-r})
\]
by Proposition \ref{vip} and Lemma \ref{unl}, along with the fact that knots are invariant under conjugation of their underlying braid.\\

Next we assume the statement is true for some $r$ and consider $w_{k-r-1}$. Without loss of generality assume $w_{k-r}=a\sigma_i^{-1} b$, and $w_{k-(r+1)}=a \sigma_i b$. We notice that our skein relation gives us a relationship between successive terms in this sequence, along with a term of lesser length. Because we know the theorem is true for shorter words as well as previous terms in the sequence, we may conclude the theorem is true always. To be precise one computes,
 
 \begin{description}
\item  $Tr_n^2(w_{k-(r+1)})=Tr_n^2(a \sigma_i b)=Tr_n^2(ba \sigma_i)$
\item  $=z \cdot Tr_n^2(ba)+Tr_n^2(ba\sigma_i^{-1}) $
\item  $=z \cdot \Delta(\hat{ba})+\Delta(\hat{ba \sigma_i^{-1}})$
\item  $=\Delta(\hat{ba \sigma_i})=\Delta(\hat{a\sigma_ib})$
\item  $=\Delta(\hat{w}_{k-(r+1)})$
 \end{description}
 \end{proof} 
 Thus, $Tr_n^2(w)=\Delta(\hat{w})$ and our proof is completed.
 \end{proof}

\section{Colored braids and linking numbers}

We conclude our paper in this section by discussing the first homomorphism family. As was alluded to previously, this family has very different properties than the other two. One will find, for example, that $\V_1$ does not satisfy a HOMFLYPT polynomial. Because of this, we will need to describe representations through $\V_1$ using alternative methods. In particular we will use coloring methods along with a representation of $\CP$ first described in \cite{FHH}.\\

\subsection{Representations through $\V_1$}
In order to discuss these representations, we must first introduce some notation.\\

We begin by creating a vector space, $V^n$. We say,
\[
 V^n = \C \text{-}span\{\mathbf{v}_S |\quad S\subseteq\{1,\ldots,n\}\}. 
\]

Furthermore, we consider the following subspaces of this vector space, 

\[ V^n_k = \C \text{-}span\{\mathbf{v}_S|\quad S\subseteq \{1,\ldots,n\} \text{ and } |S| = k\}. \]

With these defined, we introduce the following notational tools:

\begin{definition}
For a planar rook diagram $d$, let $\beta(d)$ and $\tau(d)$ denote the vertices in the top and bottom rows of $d$, respectively, which are incident to edges. We further let $f: \beta(d) \rightarrow \tau(d)$ be the function that sends a bottom vertex to its neighbor in $d$.
\end{definition}

Using these tools we may introduce a family of representations, $\rho_k: \C P_n \rightarrow End(V^n_k)$ as follows:
\[
\rho_k(d)(\mathbf{v}_{S}) = \left\{ \begin{array}{rl}
 \mathbf{v}_{f(S)} &\mbox{ if $S\subseteq \beta(d)$} \\
  0 &\mbox{ otherwise}
       \end{array} \right.
\]

It turns out that these representations are the fundamental irreducible representations of $\CP$. To be precise one has the following theorem found in \cite{FHH}:

\begin{theorem}\label{mat}
$\displaystyle\bigoplus^n_{k=0}\binom{n}{k} \rho_k$ is an isomorphism from $\C P_n$ to a direct sum of matrix algebras.
\end{theorem}

Using this theorem we will classify all representations of $\B$ through $\CP$ by using colored braids. To be precise, for any $x \in \B$ and $S \subseteq \{1,\ldots,n\}$, color each strand of $x$ green if its starting vertex is in $S$ and red otherwise. Let $r(x)$ be the number of crossings between red strands counted with sign and $r'(x)$ be the number of crossings between a red and a green strand counted with sign. We also say that for any $S \subseteq \{1,\ldots,n\}, \sum S = \displaystyle\sum_{i \in S}i$. \\

With all this defined we may proceed with the following theorem:

\begin{theorem} \label{rep}
Let $x \in B_n$ and $\pi(x) \in S_n$ be the underlying permutation of $x$, then 
\[
\rho_k \V_1(\beta)(\mathbf{v}_S)=\lambda_S^T(x)\mathbf{v}_T
\]
where $T=\pi(x)(S)$ and 
\[
\lambda_S^T(x)=(a+c+d-1)^{w(x)}(\sqrt{cd})^{w'(x)}(\sqrt{\frac{c}{d}})^{\sum T - \sum S}
\]
\end{theorem}
\begin{proof}
One notices that by the nature of group actions it suffices to show that the statement is true for $\sigma_i^{\pm1}$ where $i\in \{1,\ldots,n-1\}$.\\

Fix $S \subseteq \{1,\ldots,n\}$ and say $T = (i,i+1)S$. One verifies through quick computation that 
\begin{enumerate}
\item
\[
\rho_k \V_1(\sigma_i)(\mathbf{v}_S)=\lambda_S^T \mathbf{v}_T
\]
where
\[
\lambda_S^T(x) = \left\{ \begin{array}{rl}
  1 = (a+c+d-1)^0\sqrt{cd}^0\sqrt{\frac{c}{d}}^0 &\mbox{ if $i, i+1 \in S$}\\
  a+c+d-1 = (a+c+d-1)^1\sqrt{cd}^0\sqrt{\frac{c}{d}}^0 &\mbox{ if $i, i+1 \notin S$} \\
  c = (a+c+d-1)^0\sqrt{cd}^1\sqrt{\frac{c}{d}}^1 &\mbox{ if $i\in S, i+1 \notin S$}\\
  d = (a+c+d-1)^0\sqrt{cd}^1\sqrt{\frac{c}{d}}^{-1} &\mbox{ if $i\notin S, i+1 \in S$}
       \end{array} \right.
\]
\item
\[
\rho_k \V_1(\sigma_i^{-1})(\mathbf{v}_S)=\lambda_S^T \mathbf{v}_T
\]
where
\[
\lambda_S^T(x) = \left\{ \begin{array}{rl}
  1 = (a+c+d-1)^0\sqrt{cd}^0\sqrt{\frac{c}{d}}^0 &\mbox{ if $i, i+1 \in S$}\\
  \frac{1}{a+c+d-1} = (a+c+d-1)^{-1}\sqrt{cd}^0\sqrt{\frac{c}{d}}^0 &\mbox{ if $i, i+1 \notin S$} \\
  \frac{1}{d} = (a+c+d-1)^0\sqrt{cd}^{-1}\sqrt{\frac{c}{d}}^1 &\mbox{ if $i\in S, i+1 \notin S$}\\
  \frac{1}{c} = (a+c+d-1)^0\sqrt{cd}^{-1}\sqrt{\frac{c}{d}}^{-1} &\mbox{ if $i\notin S, i+1 \in S$}
       \end{array} \right.
\]
\end{enumerate}

As required.
\end{proof}

For the purpose of completion, one may ask whether this representation of $\B$ is related to any well known ones. Computation suggests that $\rho_k\V_1$ is some generalization of the representation discussed in \cite{FOR}; however we do not have a precise description of the connection between them.\\

\subsection{A linking number invariant}
Now that we have established the representations through $\V_1$ we turn our attentions to knot invariants through $\V_1$. We find in this case that given any braid $x \in \B$, any Markov trace will only be dependent on the linking numbers of $\hat{x}$. Despite the anti-climax of the result, it is non-trival to show as we will see.\\

We begin our proof with the following lemma:

\begin{lemma}\label{tdf}
The only trace functions $tr:\CP \rightarrow \C$ are linear combinations of matrix traces of the $V_k^n$ representations
\end{lemma}
\begin{proof}
By Theorem \ref{mat}, We know that the $\rho_k$ representations give us an isomorphism between $\CP$ and an algebra of block diagonal matrices. Because any trace functions in such an algebra would be a linear combination of traces from each block, we have our result. 
\end{proof}

This very useful fact will allow us to exhaustively consider all possible trace functions. For example, using Lemma \ref{tdf} along with Theorem \ref{rep}, we know that for any braid $x \in \B$, we only need to consider those $S \subseteq \{1,\ldots,n\}$ with $\pi(x)(S) = S$. Consider then the coloring of a braid, $x$, whose underlying permutation fixes a set $S$. We know that if a vertex was associated with a green (or red) vertex at the bottom then it will be so at the top. Furthermore, one may realize that the collection of green (or red) strands will form one or more complete components in the closure, $\hat{x}$. Finally, the scalar $\lambda^S_S(x)$ associated with the action $\rho_k\V_1(x)\mathbf{v}_S$ is determined by the total writhe of the red components of $\hat{x}$ and the total linking number between red and green components. This of course follows from Theorem \ref{rep} along with the remarks just made. Using this fact in conjunction with Lemma \ref{tdf} we conclude the following lemma:

\begin{lemma}\label{ttra}
Let $tr: \CP \rightarrow \C$ be a trace function, then for any braid $x\in \B$, $tr(\V_1(x))$ is determined by the scalars $\lambda^S_S(x)$ in 
\[
\rho_k\V_1(x)\mathbf{v}_S = \lambda^S_S(x)\mathbf{v}_S.
\]
In particular, $tr$ is at most dependent on
\begin{enumerate}
\item Which sets of strands give rise to components of $\hat{x}$,
\item the writhes of these components and
\item the linking numbers between components of $\hat{x}$.
\end{enumerate}
\end{lemma}

using this lemma we are ready to prove the main theorem of this section. We have,

\begin{theorem}
Let $Tr:\CP \rightarrow \C$ be a Markov trace function, then for any braid $x \in \B$, $Tr(\V_1(x))$ is only dependent on the linking numbers between components of $\hat{x}$.
\end{theorem}

\begin{proof}
Let $L$ be an arbitrary link of $c$ components and take $x \in \B$  such that $L = \hat{x}$. We recall that this relationship is unchanged by applying Markov moves to $x$, and thus without loss of generality we may assume that each component of $L$ uses $N$ strands of $x$ where $N$ is some sufficiently large number. This is done by conjugating $x$ to move the strands in any component of $\hat{x}$ to the right of the braid, and then appending $\sigma_n$ or $\sigma_n^{-1}$. Furthermore, each one of these moves which adds a strand to $x$ can happen using either a negative or positive crossing. This then allows us to control the writhe of each component in $\hat{x}$. One further notes that the writhe of each component will necessarily be of the opposite integer parity to the number of strands in this component. One then makes sure that through use of Markov moves, each component in $\hat{x}$ has a write of zero and $2N+1$ strands for sufficiently large $N$. From Lemma \ref{ttra} we know that $Tr$ can depend on at most the number of strands in each component, the writhes of these components and the linking number. Because $Tr$ is a Markov trace, we are free to apply the aforementioned Markov moves to any braid without changing its value. This concludes the proof.
\end{proof}

\end{document}